\documentclass[12pt]{article}
\usepackage{amsfonts}
\usepackage{amsmath}
\usepackage{amssymb}
\usepackage{amsthm}
\usepackage[english]{babel}
\begin{document}
\newtheorem{theorem}{Theorem}
\newtheorem{remark}{Remark}
\newtheorem{corollary}{Corollary}

\righthyphenmin=2
\title{On geometry of p-adic polynomials}  
\author{Evgeny Zelenov}
\date{February, 2022}
\maketitle

\begin{abstract}
An analogue of the Gauss-Lucas theorem for polynomials over the algebraic closure $\mathbb C_p$ of the field of $p$-adic numbers is considered.
\end{abstract}

\label{sec:intro}
\section{Introduction}
The pioneering works of I.~V.~Volovich \cite{1,2,3} gave impetus to the rapid development of applications of non-Archimedean analysis to models and problems of mathematical physics. The current state and bibliography can be found, for example, in the review \cite{rew}.

In this paper, we consider the geometry of non-Archimedean polynomials. The result can be helpful in the study of polynomial dynamical systems over the field of $p$-adic numbers \cite{ds1,ds2}.

The Gauss-Lucas theorem states the following. Let $P(z)\in\mathbb C[z]$ is a polynomial over the field $\mathbb C$ of complex numbers, and $P'(z)$ is its derivative. Then, all the roots of the polynomial $P'(z)$ (that is, the critical points of the polynomial $P(z)$) lie in the convex hull of the set of roots of the polynomial $P(z)$. This statement admits the following equivalent formulation. Any disk in the complex plane containing zeros of the polynomial $P(z)$ also contains all zeros of the derivative $P'(z)$ \cite{Marden}.

The purpose of this paper is to formulate and prove an analogue of this result for an algebraically closed non-Archimedean field. As such a field, we will consider the field $\mathbb C_p$ --- completion of the algebraic closure of the field $\mathbb Q_p$ of $p$-adic numbers. The norm in $\mathbb C_p$ will be denoted by $|\cdot|$.

The geometry of polynomials over non-Archimedean fields is a little-studied area. The author knows only one work on this topic \cite{Choi-Lee}, dedicated to the non-Archimedean analogue of the Sendov conjecture.

\label{sec: theorem}
\section{The main theorem}
Consider the polynomial $P(z)\in\mathbb C_p[z]$ of degree $n$ over the field $\mathbb C_p$,
\[
P(z) = a_nz^n+a_{n-1}z^{n-1}+\cdots+a_1z+a_0,\,\,a_j\in\mathbb C_p,\,\,j=0,1,\dots, n.
\]
Using $\lambda_1, \dots\lambda_n$ we denote its roots,
\[
P(z)= a_n(z-\lambda_1)(z-\lambda_2)\cdots(z-\lambda_n).
\]
Using $\omega_1, \omega_2, \dots, \omega_{n-1}$ we denote the roots of the derivative $P'(z)$,
\[
P'(z)= na_n(z-\omega_1)(z-\omega_2)\cdots(z-\omega_{n-1}).
\]
A disk of radius $r$ centered at $a\in\mathbb C_p$ is denoted by $D(a,r)$,  \[D(a,r) = \{z\in\mathbb C_p\colon |z-a|\leq r\}.\]

The following theorem is valid.
\begin{theorem}
\label{th}
Let the roots $\lambda_1,\lambda_2,\dots,\lambda_n$
of the polynomial $P(z)\in\mathbb C_p[z]$ lie within the disk $D(a,r)$. Then  each disk $D(a,r_k),\,k=1,2,\dots,n-1$, where
\[
r_k = \max\left\{r\left|\frac{j}{n}\right|^{\frac{1}{n-j}},\,j=1,2,\dots,k\right\}
\]
contains at least $k$ roots of the derivative $P'(z)$.
\end{theorem}

\begin{corollary}
\label{c1}
The disk $D\left(a,r|n|^{-1/(n-1)}\right)$ contains at least one critical point of the polynomial $P(z)$.
\end{corollary}

\begin{remark}
This statement was proved by D.~Choi and S.~Lee \cite{Choi-Lee}. Thus, the Theorem \ref{th} significantly strengthens their result.
\end{remark}

\begin{corollary}
\label{c2}
The disk $D\left(a,r|n|^{-1}\right)$ contains all critical points of the polynomial $P(z)$.
\end{corollary}

The statement follows directly from the statement of the Theorem \ref{th}  and the apparent inequality $|j|\leq1, \,j\in\mathbb Z.$

The following simple example shows that the estimate from Corollary \ref{c2} cannot be improved. Let $p=3$, consider the polynomial $P(z)=z^2(z-1)$. Its roots $\lambda_1=\lambda_2=0, \lambda_3=1$ lie within the circle $D(0,1)$. The derivative $P'(z) = 3z(z-2/3)$ has roots $\omega_1=0, \,\omega_2=2/3$ and, since $|2/3|=3$, then the minimal circle containing the roots of the derivative of the polynomial $P(z)$ is $D(0,3)$.

\begin{corollary}
\label{c3}
In the case when the order $n$ of the polynomial is not divisible by $p$ (that is, $|n|=1$), the critical points of the polynomial $P(z)$ lie within any disk containing zeros of this polynomial.
\end{corollary}
In other words, for polynomials of the order $n,\,p\nmid n$, the exact analogue of the Gauss-Lucas theorem for the field of complex numbers is valid.

\label{sec:proof}
\section{Proof of the Theorem}
Before proceeding to the proof of the Theorem \ref{th}, we will make a few remarks. First, without limiting the generality, we can put $a_n=1$ since multiplication by a nonzero constant does not change the roots of the polynomial and its derivative. We can put $a=0$ since the shift of the roots by $a$ leads to a shift of the derivative's roots by $a$. Everywhere further, we will number the roots of the polynomial $P(z)$ and the roots of its derivative $P'(z)$ in the order of non-decreasing norm:
\begin{equation}
\label{ord}
|\lambda_1|\leq |\lambda_2|\leq\cdots\leq |\lambda_n|,\,\,\,|\omega_1|\leq|\omega_2|\leq\cdots\leq|\omega_{n-1}|.
\end{equation}
Taking into account the comments made, the statement of the Theorem \ref{th} is equivalent to the following estimate for the roots of the derivative $P'(z)$,\,\, $k=1,2,\dots,n-1$:
\begin{equation}
\label{equiv}
|\omega_k|\leq r_k=r\max\left\{\left|\frac{j}{n}\right|^{\frac{1}{n-j}}, \,\,j=1,2,\dots,k\right\}.
\end{equation}
In the proof, we will use the following relations between the roots of the polynomial $P(z)$ and its coefficients (Vieta's formulas):
\begin{equation}
\label{Viete}
(-1)^{n-k}a_k = \sum_{1\leq i_1 < i_2 < \dots < i_{n-k}\leq n}\lambda_{i_1}\lambda_{i_2}\cdots\lambda_{i_{n-k}},\,\,\,k=0,1,\dots,n-1,
\end{equation}
and similar formulas for the derivative $P'(z)$:

\begin{equation}
\label{Viete1}
(-1)^{n-k}\frac{k}{n}a_k = \sum_{1\leq i_1 < i_2 < \dots < i_{n-k}\leq n-1}\omega_{i_1}\omega_{i_2}\cdots\omega_{i_{n-k}},\,\,\,k=1,\dots,n-1.
\end{equation}

The proof will be carried out by induction. First, we prove the statement of the Theorem \ref{th} for $k=1$. Indeed, a chain of relations follows from the formulas (\ref{ord}), (\ref{Viete1}), (\ref{Viete}):

\begin{multline*}
|\omega_1|^{n-1}\leq|\omega_1\omega_2\cdots\omega_{n-1}| = \\ =\frac{1}{|n|}|a_1|=\frac{1}{|n|}\left|\sum_{1\leq i_1 < i_2 < \cdots < i_{n-1}\leq n}\lambda_{i_1}\lambda_{i_2}\cdots\lambda_{i_{n-1}}\right|\leq\\\leq \frac{1}{|n|}\max_{1\leq i_1 < i_2 < \cdots < i_{n-1}\leq n}\left|\lambda_{i_1}\lambda_{i_2}\cdots\lambda_{i_{n-1}}\right|\leq\frac{1}{|n|}r^{n-1}.
\end{multline*}
In the second last inequality in the chain, we used the strong triangle inequality for the norm $|\cdot|$.
Thus, for the case of $k=1$, the Theorem \ref{th} is proved.

Now let's assume that the theorem's statement holds for $k=m-1\leq n-2$ and prove in this assumption that the theorem  holds for $k=m$ as well. Similarly to the reasoning for the case $k=1$, using the relations (\ref{Viete}), (\ref{Viete1}) and the strong triangle inequality, we obtain the following estimate for the sum of all products of $n-m$ of the roots of the polynomial $P'(z)$:

\begin{multline}
\label{est}
\left|\sum_{1\leq i_1 < i_2 < \cdots < i_{n-m}\leq n-1}\omega_{i_1}\omega_{i_2}\cdots\omega_{i_{n-m}}\right|=\left|\frac{m}{n}\right||a_{m}|= \\ =\left|\frac{m}{n}\right|\left|\sum_{1\leq i_1 < i_2 < \cdots < i_{n-m}\leq n}\lambda_{i_1}\lambda_{i_2}\cdots\lambda_{i_{n-m}}\right| \leq \left|\frac{m}{n}\right|r^{n-m}.
\end{multline}

Further proof will be carried out by contradiction. That is, suppose that the statement of the Theorem \ref{th} is not true for $k=m$. This means that the inequality (\ref{equiv}) holds for all $k\leq m-1$, but does not hold for $k=m$. Therefore, the inequality is valid:
\begin{equation}
|\omega_m|> r_m=r\max\left\{\left|\frac{j}{n}\right|^{\frac{1}{n-j}}, \,\,j=1,2,\dots,m\right\}.
\end{equation}
Under this assumption, the following chain of inequalities holds:
\begin{multline}
\label{ord1}
|\omega_1|\leq |\omega_2|\leq\cdots\leq|\omega_{m-1}|\leq r_{m-1} \leq \\ \leq r_m <|\omega_m|\leq|\omega_{m+1}|\leq\cdots\leq |\omega_{n-1}|.
\end{multline}

Therefore, among all products of $n-m$ roots of the polynomial $P'(z)$, the product $\omega_m\omega_{m+1}\cdots\omega_{n-1}$ has a strictly maximal norm.

Further, using the following property of the non-Archimedean norm:
$|a+b|=|b|$ if $|a|<|b|$, $a,b\in\mathbb C_p$, we get the estimate
\begin{equation}
\label{est1}
\left|\sum_{1\leq i_1 < i_2 < \cdots < i_{n-m}\leq n-1}\omega_{i_1}\omega_{i_2}\cdots\omega_{i_{n-m}}\right| = \left|\omega_m\omega_{m+1}\cdots\omega_{n-1}\right|\geq \left|\omega_m\right|^{n-m}.
\end{equation}
The inequalities (\ref{est}) and (\ref{est1}) directly imply the validity of the following estimate:
\[
\left|\omega_m\right|\leq r\left|\frac{m}{n}\right|^{\frac{1}{n-m}}.
\]
The last inequality contradicts assumption (2). The resulting contradiction completes the proof of the Theorem \ref{th}.

\end{document}